\documentclass{amsart}
\usepackage{amsmath,amssymb,enumerate,bbm}

\newcommand{\assign}{:=}
\newcommand{\emdash}{---}
\newcommand{\tmem}[1]{{\em #1\/}}
\newcommand{\tmop}[1]{\ensuremath{\operatorname{#1}}}
\newcommand{\tmtextbf}[1]{{\bfseries{#1}}}
\newcommand{\tmtextit}[1]{{\itshape{#1}}}
\newcommand{\tmtextup}[1]{{\upshape{#1}}}
\newenvironment{descriptioncompact}{\begin{description} }{\end{description}}
\newenvironment{enumeratenumeric}{\begin{enumerate}[1.] }{\end{enumerate}}
\newtheorem{corollary}{Corollary}
\newtheorem{definition}{Definition}
\newtheorem{example}{Example}
\newtheorem{lemma}{Lemma}
\newtheorem{proposition}{Proposition}
\newtheorem{theorem}{Theorem}

\begin{document}

\title[Formal topology and constructive mathematics]{Formal topology and
constructive mathematics: the Gelfand and
Stone-Yosida representation theorems}
\thanks{This is an expanded version of our
paper~[CS05a]{\nocite{Coquand:jucs_11_12:formal_topotoly_and_constructive}}.
For the convenience of the reader we have included more details and added a
few clarifications. There are no new results. We are grateful to Bob Lubarsky
and Fred Richman for suggesting improvements in the
presentation.}
\author{Thierry Coquand}
\email{coquand@cs.chalmers.se}
\address{Chalmers University, \\Sweden}
\author{Bas Spitters}
\thanks{Bas Spitters was
supported by the Netherlands
Organization for Scientific Research (NWO)}
\email{spitters@cs.ru.nl}
\address{Radboud University Nijmegen, \\
the Netherlands}

\begin{abstract}
  We present a constructive proof of the Stone-Yosida representation theorem
  for Riesz spaces motivated by considerations from formal topology. This
  theorem is used to derive a representation theorem for f-algebras. In turn,
  this theorem implies the Gelfand representation theorem for C*-algebras of
  operators on Hilbert spaces as formulated by Bishop and Bridges. Our proof
  is shorter, clearer, and we avoid the use of approximate eigenvalues.
\end{abstract}
\maketitle

\section{Introduction}

This paper illustrates the relevance of locale theory for constructive
mathematics. We present a constructive proof of the Stone-Yosida
representation theorem for Riesz spaces motivated by considerations from
formal topology. This theorem is used to derive a representation theorem for
f-algebras. In turn, this theorem implies the Gelfand representation theorem
for C*-algebras of operators on Hilbert spaces as formulated by Bishop and
Bridges~{\cite{Bishop/Bridges:1985}}. Our proof is shorter, \ clearer, and we
avoid the use of approximate eigenvalues.

The article is organized as follows. After dealing with some preliminaries we
prove a pointfree Stone-Yosida representation theorem for Riesz spaces. In the
next section this is used to obtain a representation theorem for f-algebras,
which in turn is used to prove the Gelfand representation theorem. Next we
discuss the similarity between Bishop's notion of compactness, i.e.~complete
and totally bounded, and compact overt spaces in formal topology. Finally we
show that the axiom of dependent choice is needed to construct points in the
formal spectrum.

We would like to stress that the present theory needs few foundational
commitments, we work within Bishop-style mathematics. Moreover, the
mathematics is predicative, even finitary, and we will not use the axiom of
choice, even countable choice, unless explicitly stated.

\section{Riesz spaces}

We present a Stone-Yosida representation theorem for Riesz spaces. This
theorem states that a Riesz space with a strong unit may be represented as a
Riesz space of functions. In fact, one can even start with an l-group, or
lattice ordered group, with a strong unit and construct a Riesz space from
this, see~{\cite{Coquand:obs}}(sec. 3).

\subsection{General definitions}

\begin{definition}
  A {\tmem{Riesz space}} $R$ is a $\mathbbm{Q}$ vector space with a compatible
  binary sup operation {\emdash} that is, such that $a + (b \vee c) = (a + b)
  \vee (a + c)$ and, moreover, $\lambda a \geqslant 0$ and $- a \leqslant 0$,
  whenever $a$ in $R$, $\lambda$ in $\mathbbm{Q}$, $a \geqslant 0$ and
  $\lambda \geqslant 0$. One can prove that A Riesz space (or {\tmem{vector
  lattice}}) is a partially ordered linear space which is a distributive
  lattice.
\end{definition}

As usual we define $a^+ \assign a \vee 0$, $a^- \assign (- a) \vee 0$, $|a|
\assign a^+ + a^-$ and $a \wedge b \assign - (- a \vee - b)$. One can prove
that a Riesz space is a lattice and that the lattice operations are compatible
with the vector space operations,
see~{\cite{Birkhoff}}{\cite{LuxemburgZaanen}}{\cite{Bourbaki}}.

\begin{definition}
  A {\tmem{strong unit 1}} in an ordered vector space $R$ is a
  positive{\footnote{We follow the standard terminology using `positive' to
  mean non-negative.}} element such that for all $a \in R$ there exists a
  natural number $n$ such that $a \leq n \cdot 1$.
\end{definition}

We will now consider a Riesz space $R$ with a strong unit. When $q$ is a
rational number, we will often write $a \leqslant q$ to mean $a \leqslant q
\cdot 1$.

\begin{definition}
  A {\tmem{representation}} of $R$ is a linear map $\sigma : R \rightarrow
  \mathbbm{R}$ such that $\sigma (1) = 1$ and $\sigma (a \vee b) = \sigma (a)
  \vee \sigma (b)$.
\end{definition}

Such a representation automatically preserves all the Riesz space structure.

\begin{example}
  If $X$ is a compact space, then $C (X)$, its space of continuous functions,
  is a Riesz space where the supremum is taken pointwise. Each point of $x$
  defines a representation $\sigma_x (f) \assign f (x) .$
\end{example}

In Example~\ref{ex:alg-of-ops} we show that a complete commutative algebra of
Hermitian operators on a Hilbert space is a Riesz space.

\subsection{Spectrum}

Given a Riesz space $R$, we will define a lattice that may be used to define
the spectrum of $R$ as a formal space, the points of which are then precisely
the representations of $R$.

Let $P$ denote the set of positive elements of a Riesz space $R$. For $a, b$
in $P$ we define $a \preccurlyeq b$ to mean that there exists $n$ such that $a
\leq nb$. We write $a \approx b$ for $a \preccurlyeq b$ and $a \succcurlyeq
b$. The following proposition is proved in~{\cite{Coquand:obs}} and involves
only elementary considerations on Riesz spaces.

\begin{proposition}
  \label{nine}We write $L (R)$ for the quotient of $P$ by $\approx$. Then $L
  (R)$ is a distributive lattice. In fact, if we define $D : R \rightarrow L
  (R)$ by $D (a) \assign [a^+]$, then $L (R)$ is the free lattice generated by
  $\{D (a) | a \in R\}$ subject to the following relations:
  \begin{enumeratenumeric}
    \item $D (a) = 0$, if $a \leq 0$;

    \item $D (1) = 1$;

    \item $D (a) \wedge D (- a) = 0$;

    \item $D (a + b) \leqslant D (a) \vee D (b)$;

    \item $D (a \vee b) = D (a) \vee D (b)$.
  \end{enumeratenumeric}
\end{proposition}

We have $D (a) \leqslant D (b)$ if and only if $a^+ \preccurlyeq b^+$ and $D
(a) = 0$ if and only if $a \leqslant 0$.

For $a$ in $R$ and rational numbers $p, q$, we write $a \in (p, q) \assign (a
- p) \wedge (q - a)$. Notice that this is an element of $R$ by our convention
that $(a - q)$ means $(a - q \cdot 1)$.

\begin{lemma}
  \label{lem:finite cover}$a \in (p, q) \preccurlyeq a \in (p, s) \vee a \in
  (t, q)$, whenever $t, s$ are rational numbers such that $t < s$. Since our
  Riesz space has a strong unit, for each $a$, there exists $p$ and $q$ such
  that $p < q$ and $a \in (p, q) = 1$. Moreover, if $I_0, \ldots, I_n$ are
  open intervals covering $(p, q)$, then $\bigvee a \in I_i \approx 1$.
\end{lemma}

\begin{proof}
  We may assume that $p < t < s < q$.
  \begin{eqnarray}
    a \in (p, s) \vee a \in (t, q) & = & ((a - p) \wedge (s - a)) \vee ((a -
    t) \wedge (q - a)) \nonumber\\
    & = & ((a - p) \vee (a - t)) \wedge ((a - t) \vee (s - a)) \wedge
    \nonumber\\
    &  & ((a - p) \vee (q - a)) \wedge ((q - a) \vee (s - a)) \nonumber\\
    & = & (a - p) \wedge ((a - t) \vee (s - a)) \wedge \nonumber\\
    &  & ((a - p) \vee (q - a)) \wedge (q - a) \nonumber\\
    & \geqslant & ((a - p) \vee (q - a)) \wedge ((a - t) \vee (s - a)) .
    \nonumber
  \end{eqnarray}
  We claim that $((a - t) \vee (s - a)) \geqslant \frac{t - s}{2}$. To show
  this we write $m \assign \frac{t + s}{2}$. Then
  \[ (a - s) \vee (t - a) = (a - m + (m - s)) \vee (t - m - (a - m)) =
     \frac{t - s}{2} + ((a - m) \vee - (a - m)) \geqslant \frac{t - s}{2} . \]
  The proof is finished by the observation that $b \preccurlyeq c$ whenever
  for some $\varepsilon$, $\varepsilon \cdot 1 \wedge b \leqslant c$.
\end{proof}

\begin{lemma}
  \label{lem:lowerbound}If $D (b_1) \vee \ldots \vee D (b_n) = 1$, then there
  exists $r > 0$ such that $D (b_1 - r) \vee \ldots \vee D (b_n - r) = 1$.
\end{lemma}

\begin{proof}
  Since $\frac{1}{N} \leqslant b_1^+ \vee \ldots \vee b_m^+$ we see that
  \begin{eqnarray*}
    \frac{1}{2 N} & \leqslant & (b_1^+ \vee \ldots \vee b_m^+) - \frac{1}{2
    N}\\
    & = & (b_1^+ - \frac{1}{2 N}) \vee \ldots \vee (b_m^+ - \frac{1}{2 N})\\
    & \leqslant & (b_1 - \frac{1}{2 N})^+ \vee \ldots \vee (b_m - \frac{1}{2
    N})^+ .
  \end{eqnarray*}
\end{proof}

The previous lemma is used to prove the following result, which can be found
as Theorem~1.11 in~{\cite{Coquand:obs}}. We will not need this, but only state
it as a motivation.

\begin{theorem}
  \label{spectrum}Define $\Sigma$, the {\tmem{spectrum}} of $R$, to be the
  locale generated by the elements $D (a)$ and the relations in
  Proposition~\ref{nine} together with the relation $D (a) = \bigvee_{r > 0} D
  (a - r)$. Then $\Sigma$ is a compact completely regular locale.
\end{theorem}

Compact completely regular locales are the pointfree analogues of compact
Hausdorff spaces. Moreover, the points of $\Sigma$ can be identified with
representations of $R$. In fact, if a representation $\sigma$ is given, then
$\sigma \in D (a)$ if and only if $\sigma (a) > 0$.

\subsection{Normable elements}

Dedekind cuts may be used to define real numbers, but it should be noted that
constructively one needs to require that such cuts $(L, U)$ are
{\tmem{located}}, i.e.~either $p \in L$ or $q \in U$, whenever $p < q$. Upper
cuts in the rational numbers may be conveniently used to deal with certain
objects that classically would also be real numbers, we call them {\tmem{upper
real}} numbers, see~{\cite{Richman:cuts}}{\cite{Vickers:locAA}}. In general,
such a cut does not have a greatest lower bound in $\mathbbm{R}$. If it does
the upper real number is called {\tmem{located}} or simply a real number.
Define the upper real $U (a) \assign \{q \in \mathbbm{Q} | \exists q' < q. a
\leq q' \cdot 1\}$ for each element of the Riesz space. If it is located $a$
is said to be {\tmem{normable}} and the greatest lower bound is denoted by
$\sup a$. Then we have $\sup a < q$ if and only if $q \in U (a)$.

\begin{proposition}
  \label{prop:normable-positive}If all elements of $R$ are normable, then the
  predicate $\tmop{Pos} (a) \assign \sup a > 0$ has the following properties:
  \begin{enumeratenumeric}
    \item If $\tmop{Pos} (a)$ and $D (a) \leqslant D (b)$, then $\tmop{Pos}
    (b)$;

    \item If $\tmop{Pos} (a \vee b)$, then $\tmop{Pos} (a)$ or $\tmop{Pos}
    (b)$;

    \item If $r$ is a strictly positive rational number, then $\tmop{Pos} (a)$
    or $D (a - r) = 0$.
  \end{enumeratenumeric}
\end{proposition}

We note that $\tmop{Pos} (a)$ if and only if $\tmop{Pos} (a^+)$.

In fact, in localic terms this shows that $\Sigma$ is {\tmem{open}}, or
{\tmem{overt}}, see~{\cite{Johnstone:open}}, but we will not need this. In
formal topology one would say that the formal space has a positivity
predicate.

We remark that the standard terminologies from the two different fields seem
to conflict. Clearly, it is not the case that $\tmop{Pos} (a)$ as soon as $a$
is positive, i.e.~$a \geqslant 0$. In particular, 0 is positive, but
$\tmop{Pos} (0)$ does not hold.

In order to prove Proposition~\ref{prop:normable-positive} we first need
three lemmas.

\begin{lemma}
  \label{lem:sup-eq}$\sup (a \vee b) = \sup a \vee \sup b$.
\end{lemma}

\begin{proof}

  \begin{descriptioncompact}
    \item[$\geqslant$] Suppose that $a \vee b \leqslant q' < q$, then $a, b
    \leqslant q'$, so both $\sup a \leqslant q'$ and $\sup b \leqslant q'$.
    Thus $q' \geqslant \sup a \vee \sup b$.

    \item[$\leqslant$] If $q' \geqslant a, b$, then $q' \geqslant a \vee b$
    and hence $\sup a \vee b \leqslant \sup a \vee \sup b$.
  \end{descriptioncompact}
\end{proof}

\begin{lemma}
  \label{lem:eleven}Let $r$ be a rational number. If $\sup b < r$ and $r <
  \sup (b \vee c)$, then $r < \sup c.$
\end{lemma}

\begin{proof}
  If $\sup b < r$, then $b \leqslant r' < r$, for some rational number $r'$.
  So, using Lemma~\ref{lem:sup-eq}
  \[ r < \sup (b \vee c) = \sup b \vee \sup c \leqslant r' \vee \sup c. \]
  Consequently, $r < \sup c$.
\end{proof}

\begin{lemma}
  \label{lem:sup}If $0 < \sup (b_1 \vee b_2)$, then $\sup b_1 > 0$ or $\sup
  b_2 > 0$.
\end{lemma}

\begin{proof}
  Suppose that $0 < r < \sup (b_1 \vee b_2)$. Either $\sup b_1 < r$ or $\sup
  b_1 > 0$. In latter case we are done. In the former case $\sup b_2 > r$ by
  Lemma~\ref{lem:eleven}.
\end{proof}

\begin{proof}
  \tmtextbf{[of Proposition~\ref{prop:normable-positive}]}

  Property 1 is clear.

  Property 2 is Lemma~\ref{lem:sup}.

  Finally, to prove property 3 we decide whether $\sup a > 0$ or $\sup a < r$.
  In the former case $\tmop{Pos} (a)$. In the latter case $D (a - r) = 0$.
\end{proof}

\begin{corollary}
  \label{cor:pos}If $D (a_1) \vee \ldots \vee D (a_n) = 1$, we can find $i_1 <
  \ldots < i_k$ such that $D (a_{i_1}) \vee \ldots \vee D (a_{i_k}) = 1$ and
  $\tmop{Pos} (a_{i_1}), \ldots, \tmop{Pos} (a_{i_k})$.
\end{corollary}

\begin{proof}
  By Lemma~\ref{lem:lowerbound} there exists $r > 0$ such that $D (a_1 - r)
  \vee \ldots \vee D (a_n - r) = 1$. By
  Proposition~\ref{prop:normable-positive} for all $i$, $D (a_i - r) = 0$ or
  $\tmop{Pos} (a_i)$. From this the result follows.
\end{proof}

\begin{lemma}
  \label{lem:intersection}Let $I : = (p, q)$, $J \assign (r, s)$. Define $I +
  J \assign (p + r, q + s)$ and $I \vee J \assign (p \vee r, q \vee s)$. If
  $|a + b - c| \leqslant \varepsilon$ and $\tmop{Pos} (a \in I \wedge b \in J
  \wedge c \in K)$, then the distance between $I + J$ and $K$ is bounded by
  $\varepsilon$. If $|a \vee b - c| \leqslant \varepsilon$ and $\tmop{Pos} (a
  \in I \wedge b \in J \wedge c \in K)$, then the distance between $I \vee J$
  and $K$ is bounded by $\varepsilon$.

  Finally, if $|a - b| \leqslant \varepsilon$ and $\tmop{Pos} (a \in I \wedge
  b \in J)$, then the distance between $I$ and $J$ is bounded by
  $\varepsilon$.
\end{lemma}

\begin{proof}
  We only prove the last fact. If the distance between $I$ and $J$ is bigger
  than $\varepsilon$, than $r - q > \varepsilon$ or $p - s > \varepsilon$.
  Consequently, $D (q - a \wedge b - r) = 0$ or $D (a - p \wedge s - b) = 0$.
  Both cases imply that $D (a \in I \wedge b \in J) = 0$.
\end{proof}

\begin{definition}
  A Riesz space $R$ is {\tmem{separable}} if there exists a sequence $a_n$
  such that for all $a$ and \ $\varepsilon > 0$, there exists $n$ such that
  $|a - a_n | \leqslant \varepsilon$.
\end{definition}

\begin{theorem}
  \label{point}\tmtextup{[DC]} Let $R$ be a separable Riesz space all elements
  of which are normable. Assume that $\tmop{Pos} (a)$, then there exists a
  representation $\sigma$ such that $\sigma (a) > 0.$
\end{theorem}

\begin{proof}
  We write $\varepsilon_n \assign 2^{- n}$. Using dependent choice and
  Lemma~\ref{lem:finite cover} we define a sequence $(q_n)$ of rationals such
  that
  \[ \tmop{Pos} (a \in ( \frac{\sup a}{2}, \sup a) \wedge a_0 \in (q_0 -
     \varepsilon_0, q_0 + \varepsilon_0) \wedge \ldots \wedge a_n \in (q_n -
     \varepsilon_n, q_n + \varepsilon_n)) . \]
  If $b \in R$, we can find a sequence of elements $a_{n_k}$ such that for any
  $\varepsilon > 0$ we have $|b - a_{n_k} | \leq \varepsilon$ when $k$ is
  large enough. Then $q_{n_k}$ is a Cauchy sequence and we define $\sigma (b)
  \assign \lim_k q_{n_k}$. By Lemma~\ref{lem:intersection} this definition
  does not depend on the choice of the sequence $a_{n_k}$. The map $\sigma$ is
  a representation such that $\sigma (a) > 0$ and $\sigma (a_n) \in (q_n -
  \varepsilon_n, q_n + \varepsilon_n)$ for all $n$.
\end{proof}

A suggestive way to state that $\sigma (a) > 0$ is to say that $\sigma$ is a
{\tmem{point}} in $D (a)$.

Let $\Sigma$ be the set of representations of $R$. We call $\Sigma$ the
{\tmem{spectrum}} of $R$. Each representation is a bounded linear functional.
Each element $a$ of $R$ defines a pseudo norm $\rho_a (\varphi) \assign |
\varphi (a) |$ on the space of bounded linear functionals. If $R$ is separable
and $a_n$ is a dense sequence in $\{a \in R : |a| \leqslant 1\}$, we can
collect, like Bishop, all the pseudo-norms into one norm $\rho (\varphi)
\assign \sum_n 2^{- n} | \varphi (a_n) |$. Considering the restriction of
these pseudo norm to the spectrum, which is not a linear space, we obtain a
pseudo metric.

We have the following Stone-Yosida representation theorem,
see~{\cite{Stone:spectrum2}}{\cite{Yosida}}.

\begin{theorem}
  \tmtextup{\label{Gelfand}[DC]} Let $R$ is a separable Riesz space all
  elements of which are normable. The spectrum $\Sigma$ is a complete totally
  bounded metric space. For $a$ in $R$, we define $\hat{a} : \Sigma
  \rightarrow \mathbbm{R}$ by $\hat{a} (\sigma) \assign \sigma (a)$ and $\|a\|
  \assign \sup (|a|)$. Then $\sup_{\sigma} | \hat{a} (\sigma) | =\|a\|.$
  Finally, the set of functions $\hat{a}$ is dense in $C (\Sigma)$.
\end{theorem}

\begin{proof}
  We define $U_{a r s} \assign D (a \in (r, s))$ as an element of the lattice
  $L (R)$. Let $\varepsilon > 0$ and $a_1, \ldots, a_n$ in $R$. For each $i$
  we construct, using Lemma~\ref{lem:finite cover}, finitely many $s_{i j},
  r_{i j}$ such that $s_{i j} - r_{i j} < \varepsilon$ and $\bigvee_j U_{a_i
  r_j s_j} = 1$ in the lattice $L (R)$. By Corollary~\ref{cor:pos} one can
  assume all these elements to be positive. Each of them contains a point
  $\sigma_{i j}$ by Theorem~\ref{point}. If $\tau$ in $U_{a_i r_j s_j}$, then
  $\rho_{a_i} (\tau, \sigma_{i j}) = | \tau (a) - \sigma_{i j} (a) | < s_{i j}
  - r_{i j}$. Since there are finitely many $U_{a_i r_j s_j}$ covering
  $\Sigma$, the collection of these points forms an $\varepsilon$-net for the
  pseudonorms $\rho_{a_i}$. Consequently, $\Sigma$ is totally bounded.

  It is straightforward to show that $\Sigma$ is also complete as a uniform
  space.

  For each $\sigma \in \Sigma$ we have $| \sigma (a) | \leqslant \|a\|$. To
  see this suppose that $\sigma (a) >\|a\|$. Then there exists $\varepsilon >
  0$ such that $\sigma (a) - a \geqslant \varepsilon 1$, however $\sigma
  (\sigma (a) - a) = 0$. If $r <\|a\|$, then by Theorem~\ref{point}, there
  exists $\sigma$ such that $r < | \sigma (a) |$.

  Finally, the density follows from Proposition~3.1 in~{\cite{Coquand:obs}}.
  Its proof involves only elementary properties of Riesz spaces.
\end{proof}

Notice the interplay between the {\tmem{pointwise}} and {\tmem{pointfree}}
framework. From a formal covering, in the lattice $L (R)$, it is possible to
deduce that $\Sigma$, a metric space, is totally bounded. This is remarkable
since there are examples of Riesz spaces $R$ such that in a recursive
interpretation of Bishop's mathematics $\Sigma$ does not have enough points.
For instance, consider the Riesz space $R$ of continuous real functions on
Cantor space ($2^{\omega}$). The spectrum of $R$ is precisely Cantor space and
the representations are its points{\footnote{To see this, note that the
characteristic function $\chi_u$ of any basic open $u$ is continuous. So,
given a representation $\sigma$, $\sigma (\chi_u \wedge (1 - \chi_u)) = 0$ and
$\sigma (\chi_u \vee (1 - \chi_u)) = 1$. Consequently, $\sigma (\chi_u)$ is
either 0 or 1. It follows that the representations can be identified with the
points.}}, which is known not to have enough points in a recursive
interpretation. In particular, this means that we have a collection of open
sets which covers all the recursive points, but does not allow a finite
subcover. This is possible since this collection does not cover the space in
the usual terminology of formal topology.

\section{f-algebras}

In this section we apply the results of the previous section to f-algebras.

\begin{definition}
  An {\tmem{f-algebra}} is a Riesz space with a strong unit and a
  commutative{\footnote{\tmtextup{One can prove classically that the
  commutativity requirement follows form the other properties of an f-algebra.
  We intend to provide a constructive proof of this separately using the
  pointfree description of the spectrum.}}} multiplication such that $0
  \leqslant ab$, whenever $0 \leqslant a$ and $0 \leqslant b$.
\end{definition}

\subsection{f-algebra of operators}

\begin{example}
  \label{ex:alg-of-ops}If $R$ is a complete commutative algebra of normable
  self-adjoint operators on a Hilbert space $H$, then $R$ is a Riesz space
  with the order $\leqslant$ defined by $0 \leqslant A$ if and only if $(A u,
  u) \geqslant 0$ for all $u$ in $H$.
\end{example}

In the rest of this subsection we prove that if $A, B \geqslant 0,$ then $AB
\geqslant 0$.

We have now defined two notions of boundedness on the algebra of operators.
One as a bounded operator: $A$ is bounded by $a$ if for all $x$, $\|A x\|^2
\leqslant a\|x\|^2$. The other from the ordering: $A$ is bounded by $a$ if $A
\leqslant a I$, where $I$ is the identity operator.

\begin{lemma}
  \label{lem:norm}The two notions of boundedness coincide {\emdash} that is,
  for all $x$, $(A x, x) \leqslant (a x, x)$ if and only if for all $x$, $\|A
  x\|^2 \leqslant a^2 \|x\|^2$. Consequently, $\|A^2 \|=\|A\|^2$.
\end{lemma}

\begin{proof}
  The usual proof, for instance in~{\cite{Lang}} using the polarization
  identity, is constructive.
\end{proof}

Since $(AB^2 x, x) = (A (B x), (B x))$, we see that $A B^2 \geqslant 0$,
whenever $A \geqslant 0$. This suffices to prove that $R$ is an ordered ring.

\begin{lemma}
  \tmtextit{\label{lem:Riesz}\tmtextup{{\cite{Riesz:1930}} (p33, footnote
  9)}} Let $R$ be a as above. Then every positive element is the uniform limit
  of a sum of squares.
\end{lemma}

\begin{proof}
  We can assume that $0 \leqslant A \leqslant 1.$ Define $A_0 : = A$ and $A_{n
  + 1} : = A_n - A_n^2$. Then $0 \leqslant A_{n + 1} \leqslant 1,$ since $A_{n
  + 1} = A_n (1 - A_n)^2 + (1 - A_n) A_n^2 \geqslant 0$ and $1 - A_{n + 1} = 1
  - A_n + A_n^2$. Moreover, $A_{n + 1} = A_n - A_n^2 \leqslant A_n$. Since $A
  = A_1^2 + \cdots + A_n^2 + A_{n + 1},$ we have $A_n^2 \leqslant 1 / n \to
  0$. By Lemma~\ref{lem:norm} this implies that $A_n \rightarrow 0$.
\end{proof}

\begin{corollary}
  $AB \geqslant 0$, whenever $A, B \geqslant 0$.
\end{corollary}

\begin{proof}
  If $A, B \geqslant 0$, then $(A B x, x) = \sum (A B_n x, B_n x) \geqslant
  0$, where $B_n$ is a sequence such that $\Sigma B_n^2$ converges to $B$.
\end{proof}

The following lemma shows that one can construct the square
root{\footnote{This is the usual lemma that $R$ admits square root of positive
elements if $R$ is complete. Notice that the proof is directly constructive,
and it corresponds to the usual Taylor expansion of $(1 - x)^{1 / 2}$.}} when
$R$ is complete and thus one can define the absolute value as $|A| \assign
\sqrt{A^2}$. From the absolute value one first defines $A^+ \assign (|A| + A)
/ 2$ and then $A \vee B \assign A + (B - A)^+$. Consequently, the algebra is a
Riesz space and an f-algebra.

\begin{lemma}
  \label{sequence}For all $A \geqslant 0$ we can build a Cauchy sequence
  $(A_n)$ of positive elements such that $A_n^2 \rightarrow A$.
\end{lemma}

\begin{proof}
  We can assume $0 \leq A \leq I$. We define the two sequences $A_n \in [0,
  I]$ and $r_n \in [0, 1]$ defined by $A_0 = 0$ and $r_0 = 0$ and
  \[ A_{n + 1} = \frac{1}{2} (1 - A + A_n^2) \hspace{0.25em} \hspace{0.25em}
     \hspace{0.25em} \hspace{0.25em} \hspace{0.25em} \hspace{0.25em}
     \hspace{0.25em} \hspace{0.25em} \hspace{0.25em} \hspace{0.25em}
     \hspace{0.25em} \hspace{0.25em} \hspace{0.25em} \hspace{0.25em}
     \hspace{0.25em} r_{n + 1} = \frac{1}{2} (1 + r_n^2) \]
  Clearly, we have $A_n \leq r_n$ for all $n$.

  We claim that we have for all $n$
  \[ A_n \leq A_{n + 1} \hspace{2em} r_n \leq r_{n + 1} \hspace{2em} A_{n + 1}
     - A_n \leq r_{n + 1} - r_n \]
  This is proved by induction from the equalities
  \[ A_{n + 1} - A_n = \frac{1}{2} (A_n + A_{n - 1}) (A_n - A_{n - 1})
     \hspace{2em} r_{n + 1} - r_n = \frac{1}{2} (r_n + r_{n - 1}) (r_n - r_{n
     - 1}) \]
  It follows that we have
  \[ (I - A_n)^2 - A = 2 (A_{n + 1} - A_n) \leq 2 (r_{n + 1} - r_n) \]
  In order to conclude, all is left is to show that $(r_n)$ has limit $1$. We
  know that $0 \leq r_n \leq r_{n + 1} \leq 1$ and we have
  \[ 1 - r_{n + 1} = \frac{1}{2} (1 - r_n^2) = (1 - r_n) \frac{1}{2} (1 + r_n)
     \leq (1 - r_n) (1 - \frac{\epsilon}{2}) \]
  if $r_n \leq 1 - \epsilon$. This shows that if $(1 - \frac{\epsilon}{2})^N
  \leq \epsilon$ we have $1 - r_n \leq \epsilon$ for all $n \geq N$.
\end{proof}

\subsection{Gelfand representation}

Any f-algebra is a Riesz space so we have a Gelfand representation of the
f-algebra qua Riesz space, see Theorem~\ref{Gelfand}.

\begin{theorem}
  The Gelfand transform $\hat{\cdot}$ preserves multiplication.
\end{theorem}

\begin{proof}
  Since $2 ab = (a + b)^2 - a^2 - b^2$. We need to prove that $\hat{\cdot}$
  preserves squares {\emdash} that is $\sigma (a^2) = \sigma (a)^2$. For this
  we first prove: $\sigma (ab) > 0$, whenever $\sigma (a), \sigma (b) > 0$.

  If $\sigma (a) \geqslant r > 0$, then $\sigma (a - r) > 0$. By Lemma~6.3
  in~{\cite{Coquand:obs}}, $(a - r)^+ \wedge b^+ \leqslant \frac{1}{r}
  (ab)^+$, so $\sigma (ab)^+ \geqslant r (\sigma (b)^+ \wedge \sigma (a -
  r)^+) > 0$, which was to be proved.

  Suppose that $| \sigma (a) | < q$. Then $q > \sigma (a)$ and so $\sigma (q -
  a) > 0$. Similarly, $\sigma (q + a) > 0$. Consequently, $\sigma (q^2 - a^2)
  = \sigma ((q - a) (q + a)) > 0$. By a similar argument we see that if $|
  \sigma (a) | < q$, then $| \sigma (a^2) | < q^2$. We conclude that $\sigma
  (a^2) = \sigma (a)^2$.
\end{proof}

We have proved the following representation theorem for f-algebras which
explains the name f-algebra: an f-algebra is an abstract function algebra.

\begin{theorem}
  \tmtextup{\label{thm:bish}[DC]}Let $\mathcal{A}$ be a separable f-algebra of
  normable elements, then the spectrum $\Sigma$ is a compact metric space and
  there exists an f-algebra embedding of $\mathcal{A}$ into $C (\Sigma)$.
\end{theorem}

We now specialize this theorem to the f-algebra in Example~\ref{ex:alg-of-ops}
and obtain Bishop's version of the Gelfand representation theorem. In fact,
like Bishop we first prove the theorem for Hermitian operators. As a corollary
we obtain the Gelfand duality theorem for a separable Abelian C*-algebra,
exactly as stated by Bishop~{\cite{Bishop/Bridges:1985}} Cor.8.28 by
considering its self-adjoint part which is an f-algebra.

\begin{corollary}
  \tmtextup{[DC]}Let $\mathcal{A}$ be a separable f-algebra of normable
  Hermitian operators on a Hilbert space, then the spectrum is a compact
  metric space and there exists an f-algebra embedding of $\mathcal{A}$ into
  $C (\Sigma)$.
\end{corollary}

\begin{theorem}
  \tmtextup{\label{thm:Gelfand}[DC]}Let $R$ be an Abelian C*-algebra of
  operators on a Hilbert space. Then there exists a C*-algebra embedding
  $\varphi$of $R$ into $C (\Sigma, \mathbbm{C})$, where $\Sigma$ is a compact
  metric space. Moreover, $\varphi (1) = 1$ and $R$ is norm-dense.
\end{theorem}

Bishop's Gelfand representation theorem states that for any commutative
algebra of normable operators on a separable Hilbert space there exists a
norm-preserving isomorphism to the algebra of continuous functions on its
spectrum. To prove that this map is norm-preserving Bishop proves that certain
$\varepsilon$ eigenvectors can be computed. In fact, the computational
information of the $\varepsilon$ eigenvectors is used only to prove the
non-computational statement that the map is norm-preserving. In contrast, we
work directly on the approximations so that we can avoid these unused
computational steps.

\subsection{Peter-Weyl}

For a typical application, we let $G$ be a compact group and $R$ be the
algebra of operators over $L_2 (G)$ generated by the unit operator and the
operators $T (f) (g) \assign f \ast g$, where $\ast$ denotes the convolution
product. Each operator $T (f)$ is compact and hence normable. The non-trivial
representations of $R$ are then exactly the characters of the group $G$. This
gives a reduction of the Peter-Weyl theorem to the Gelfand representation
theorem, see~{\cite{CoquandSpittersPW}}.

\section{\label{sub:compactovert}Compact overt locales}

Bishop defines a metric space to be compact if it is complete and totally
bounded and proves that all uniformly continuous functions defined on such a
metric space are normable. In contrast, in the framework of locale theory it
is not true in general that all functions on a compact regular locale are
normable, i.e.~the norm is only defined as an upper real which may not be a
located.

However, the locale considered in Theorem~\ref{spectrum} is not only compact
completely regular, but also {\tmem{overt}} {\emdash} that is, has a
positivity predicate {\emdash} as shown by
Proposition~\ref{prop:normable-positive}. In this case all the continuous
functions are normable. This is a general fact.

\begin{theorem}
  \label{thm:open-sup}If $X$ is a compact completely regular locale, then $X$
  is overt if and only if for any $f \in C (X)$ there exists $\sup f \in
  \mathbbm{R}$ such that $\sup f < s$ if and only if $f^{- 1} (- \infty, s) =
  X.$
\end{theorem}

\begin{proof}
  We prove only the `only if' part. If $X$ is overt and $f \in C (X)$, then an
  approximation of the supremum can be found by considering a finite covering
  of $X$ by positive opens of the form $f^{- 1} (r, s)$, where $s - r$ is
  small.
\end{proof}

The previous facts suggest a similarity between Bishop's compact metric spaces
and the compact overt spaces in formal topology.

\begin{center}
  \begin{tabular}{l}
    Bishop compact $\Leftrightarrow$ compact overt
  \end{tabular}
\end{center}

Clearly this requires further developments building on ideas
in~{\cite{MartinLof:NCM}}{\cite{Johnstone:open}}. However, we postpone this to
further work.

It is interesting to note that Paul Taylor has independently found a similar
relation between Bishop compact and compact overt in the context context of
his abstract Stone duality~{\cite{Taylor}}. He also introduced the term
\tmtextit{overt}.

We would like to conclude this discussion with the following comparison
between the three following frameworks: classical mathematics with the axiom
of choice, Bishop's mathematics and our framework, predicative constructive
mathematics without dependent choice{\footnote{This framework is related to
Richman's proposal to develop constructive mathematics without using countable
choice, see~{\cite{Richman:fta}}.}}. Using classical logic and the axiom of
choice one can show that the spectrum defined in Theorem~\ref{spectrum} has
enough points~{\cite{johnstone:stone}}. Thus in this setting the pointfree and
pointwise description of the space coincide. In a recursive interpretation of
Bishop's framework these descriptions differ. However, using dependent choice,
normability and separability assumptions, we have shown how to deduce that
$\Sigma$ is totally bounded from the pointfree description of $\Sigma$.

Our conclusion is that the best formulation of the representation theorem is
the pointfree one, since, besides being neutral on the use of the axiom of
choice and classical logic, it implies the usual formulations both in Bishop's
framework and in classical mathematics.

\section{Choice}

\subsection{No points}

As mentioned before, when all the elements of the algebra are normable, one
can construct points in the spectrum using dependent choice. We claim that
dependent choice is needed for this. In fact, it is known that there exist
compact overt locales for which we need countable choice to construct a point.
Thus it suffices to consider the space of continuous functions on such a
locale. We think that a nice example can be extracted
from~{\cite{Richman:fta}}.

Richman~{\cite{Richman:fta}}(p.5) gave an informal argument that indicates
that one can not construct the zeroes of the complex polynomial $X^2 - a$
unless one knows whether $a = 0$ or not. To $a$ in $\mathbbm{C}$ we can
associate the locale $Y_a$ of roots of $X^2 - a$. The existence of a point in
$Y_a$ requires dependent choice. On the other hand using results
from~{\cite{Vickers:locAA}} it should be possible to show that $Y_a$ is
compact overt as an element of the completion of the metric space
$n$-multisets in $\mathbbm{C}$. The metric on this space is the usual
Hausdorff metric on compact subsets of $\mathbbm{C}$.

Finally, we remark that this leaves open the question whether it is possible
to construct the points of the spectrum of a discrete countable Riesz space
over the rationals, or more generally, to construct the points of the spectrum
of a separable Riesz space.

\subsection{Spreads}

It is interesting to note that Richman~{\cite{Richman:spreads}} proposes to
use spreads to avoid dependent choice. We suggest to use formal spaces
instead. One motivation of formal
topology~{\cite{MartinLof:NCM}}{\cite{Sambin:1987}} was precisely to give a
direct treatment of Brouwer's spreads by working with trees of finite
sequences. Formal topology may be seen as a predicative and constructive
version of locale theory. Johnstone~{\cite{johnstone:stone}} stresses that one
may avoid the use of the axiom of choice in topology by using locale theory
and dealing directly with the opens. In this light it may not be so surprising
that Richman uses spreads to avoid choice.

Richman's definition of spread differs in two respects from Heyting's
definition. The branching of the tree is arbitrary, i.e.~not necessarily
indexed by the natural numbers, and it is not decidable whether or not a
branch can be continued. This may be compared to the present situation where
we study the maximal spectrum $\Sigma$. When $\Sigma$ has a countable base, we
may define it as a finitely branching tree. When furthermore $\Sigma$ is
overt, every positive branch can be continued in a positive way. One
difference between Richman's spreads and our approach is that Richman requires
the infinite branches to be elements of a metric space.

\section{Conclusion}

We gave a constructive proof of the Stone-Yosida representation theorem for
Riesz spaces. This theorem was used to prove a representation theorem for
f-algebras, from which we derived the Gelfand representation theorem for
commutative C*-algebras of operators on a Hilbert space. This constructive
theorem generalizes the one by Bishop and Bridges. In a similar way one may
prove a generalization of Bishop's spectral theorem,
see~{\cite{spitters:obs}}.

It should be noted that we have used normability and separability hypothesis
in the statements of the main theorems and used the axiom of dependent choice.
In fact, without these hypothesis we can still obtain the spectrum as a
compact locale. The normability is necessary to show that the spectrum is
overt. The separability hypothesis is used to obtain a metric space instead of
a uniform space. Finally, the axiom of dependent choice is used in
Theorem~\ref{point} to construct a point in each positive open, and thus
obtain a metric space in the sense of Bishop.

In this context, we would like to mention a problem for both constructive
versions of the Gelfand representation theorem: can it be applied to construct
the Bohr compactification of, say, the real line, like
Loomis~{\cite{Loomis:AHA}}? Considering that the Stone-�ech compactification
has been successfully treated in locale theory~{\cite{johnstone:stone}}, one
would hope that a similar treatment is possible. Since the almost periodic
functions do not form an algebra constructively, we may consider the f-algebra
of functions generated by them. However, in this algebra not all elements are
normable. Any element of the group determines a point in the spectrum.
However, it is not possible to extend the group operation to the spectrum and
obtain a localic group, since every compact localic group has a positivity
predicate~{\cite{Wraith:groups}} and since, moreover, the spectrum is compact
this would imply that all the functions in the f-algebra are normable. This,
as we stated before, is not the case. See Spitters~{\cite{Spitters:ap}} and
the references therein for a constructive theory of almost periodic functions.

\bibliographystyle{alpha}\bibliography{douglas.bib,StoneYosida.bib,expanded.bib}

\begin{thebibliography}{Spi05b}

\bibitem[BB85]{Bishop/Bridges:1985}
Errett Bishop and Douglas Bridges.
\newblock {\em {Constructive analysis}}, volume 279 of {\em {Grundlehren der
  Mathematischen Wissenschaften}}.
\newblock Springer-Verlag, 1985.

\bibitem[Bir67]{Birkhoff}
Garrett Birkhoff.
\newblock {\em Lattice theory}.
\newblock Third edition. American Mathematical Society Colloquium Publications,
  Vol. XXV. American Mathematical Society, Providence, R.I., 1967.

\bibitem[Bou64]{Bourbaki}
N.~Bourbaki.
\newblock {\em Alg\`ebre, chapitre 6.}
\newblock Hermann, 1964.

\bibitem[Coq05]{Coquand:obs}
Thierry Coquand.
\newblock About {S}tone's notion of spectrum.
\newblock {\em J. Pure Appl. Algebra}, 197(1-3):141--158, 2005.

\bibitem[CS05a]{Coquand:jucs_11_12:formal_topotoly_and_constructive}
T.~Coquand and B.~Spitters.
\newblock Formal {T}opology and {C}onstructive {M}athematics: the {G}elfand and
  {S}tone-{Y}osida {R}epresentation {T}heorems.
\newblock {\em Journal of Universal Computer Science}, 11(12):1932--1944, 2005.
\newblock
  \verb|http://www.jucs.org/jucs_11_12/formal_topotoly_and_constructive|.

\bibitem[CS05b]{CoquandSpittersPW}
Thierry Coquand and Bas Spitters.
\newblock A constructive proof of the {P}eter-{W}eyl theorem.
\newblock {\em Mathematical Logic Quarterly}, 4:351--359, 2005.

\bibitem[Joh82]{johnstone:stone}
Peter~T. Johnstone.
\newblock {\em Stone Spaces}.
\newblock Number~3 in Cambridge studies in advanced mathematics. Cambridge
  University press, 1982.

\bibitem[Joh84]{Johnstone:open}
Peter~T. Johnstone.
\newblock Open locales and exponentiation.
\newblock In {\em Mathematical applications of category theory (Denver, Col.,
  1983)}, volume~30 of {\em Contemp. Math.}, pages 84--116. Amer. Math. Soc.,
  Providence, RI, 1984.

\bibitem[Lan83]{Lang}
Serge Lang.
\newblock {\em Real analysis}.
\newblock Addison-Wesley Publishing Company Advanced Book Program, Reading, MA,
  second edition, 1983.

\bibitem[Loo53]{Loomis:AHA}
Lynn~H. Loomis.
\newblock {\em An introduction to Abstract Harmonic Analysis}.
\newblock University Series in Higher Mathematics. van Nostrand, New York,
  1953.

\bibitem[LZ71]{LuxemburgZaanen}
W.~A.~J. Luxemburg and A.~C. Zaanen.
\newblock {\em Riesz spaces. {V}ol. {I}}.
\newblock North-Holland Publishing Co., Amsterdam, 1971.
\newblock North-Holland Mathematical Library.

\bibitem[ML70]{MartinLof:NCM}
Per Martin-L\"of.
\newblock {\em Notes on constructive mathematics}.
\newblock Almqvist \& Wiksell, Stockholm, 1970.

\bibitem[Ric98]{Richman:cuts}
Fred Richman.
\newblock Generalized real numbers in constructive mathematics.
\newblock {\em Indagationes Mathematicae}, 9:595--606, 1998.

\bibitem[Ric00]{Richman:fta}
Fred Richman.
\newblock The fundamental theorem of algebra: a constructive development
  without choice.
\newblock {\em Pacific Journal of Mathematics}, 196:213--230, 2000.

\bibitem[Ric02]{Richman:spreads}
Fred Richman.
\newblock Spreads and choice in constructive mathematics.
\newblock {\em Indagationes Mathematicae}, 13:259--267, 2002.

\bibitem[Rie32]{Riesz:1930}
F.~Riesz.
\newblock Ueber die linearen {T}ransformationen des komplexen {H}ilbertschen
  {R}aumes.
\newblock {\em Acta Sci. Math.}, 5:23--54, 1930-32.

\bibitem[Sam87]{Sambin:1987}
Giovanni Sambin.
\newblock Intuitionistic formal spaces - a first communication.
\newblock In D.~Skordev, editor, {\em Mathematical logic and its Applications},
  pages 187--204. Plenum, 1987.

\bibitem[Spi05a]{Spitters:ap}
Bas Spitters.
\newblock Almost periodic functions, constructively.
\newblock {\em Logical Methods in Computer Science}, 1(3:3):1--7, 2005.

\bibitem[Spi05b]{spitters:obs}
Bas Spitters.
\newblock Constructive algebraic integration theory without choice.
\newblock {\em Dagstuhl proceedings}, 2005.

\bibitem[Sto41]{Stone:spectrum2}
M.~H. Stone.
\newblock A general theory of spectra. {II}.
\newblock {\em Proc. Nat. Acad. Sci. U. S. A.}, 27:83--87, 1941.

\bibitem[Tay05]{Taylor}
Paul Taylor.
\newblock A lambda calculus for real analysis.
\newblock In Tanja Grubba, Peter Hertling, Hideki Tsuiki, and Klaus Weihrauch,
  editors, {\em CCA}, volume 326-7/2005 of {\em Informatik Berichte}, pages
  227--266. FernUniversit{\"a}t Hagen, Germany, 2005.

\bibitem[Vic05]{Vickers:locAA}
Steven Vickers.
\newblock Localic completion of generalized metric spaces {I}.
\newblock {\em Theory and Applications of Categories}, 14:328--356, 2005.

\bibitem[Wra90]{Wraith:groups}
G.~C. Wraith.
\newblock Unsurprising results on localic groups.
\newblock {\em J. Pure Appl. Algebra}, 67(1):95--100, 1990.

\bibitem[Yos42]{Yosida}
K{\^o}saku Yosida.
\newblock On the representation of the vector lattice.
\newblock {\em Proc. Imp. Acad. Tokyo}, 18:339--342, 1942.

\end{thebibliography}

\end{document}